\documentclass{article}[12pt]

\newtheorem{thm}{Theorem}
\newtheorem{cor}{Corollary}

\newtheorem{pro}{Proposition}
\newtheorem{rem}{Remark}

\usepackage[ marginparsep=5mm,heightrounded,centering,margin=3cm]{geometry}

\newcommand{\di}{\displaystyle}

\usepackage{enumerate}
\usepackage{amssymb,amsmath}

\usepackage[numbers]{natbib}
\newenvironment{proof}[1][Proof:]{\begin{trivlist}
\item[\hskip \labelsep {\bfseries #1}]}{\end{trivlist}}
\usepackage{geometry}
\title{Finitely presented quadratic algebras of intermediate growth}
\author{Dilber Ko\c{c}ak}
\date{06.03.2015}
\begin {document}
\maketitle

\begin{abstract} In this article, we give two examples of finitely presented 
 quadratic algebras (algebras presented by quadratic relations) of intermediate growth.
 \end{abstract}

\section{Introduction}

Let $A$ be a finitely generated algebra over a field $k$
with generating set $S=\{x_1,\ldots, x_m\}$
.
We denote by 
$A_n$ the subspace of elements of degree at most $n$, then $A=\bigcup_{n=0}^{\infty}A_n$. 
The growth function $\gamma_{A}^{S}$
of $A$ with respect to $S$ is defined as the dimension of the vector space
$A_n$ over $k$, 
$$\gamma_{A}^{S}(n)=dim_{k}(A_n)$$
The function $\gamma_{A}^{S}$ depends on the generating set $S$. 
This dependence can be removed
by introducing an equivalence relation: Let $f$ and $g$ be eventually monotone increasing 
and positive 
valued functions on $\mathbb{N}$. Set $f\preceq g$ if and only if there exist $N>0$,
$C>0$, such that 
$f(n)\leq g(Cn)$, for $n\geq N$, and $f \sim g $ if and only if $f \preceq g$ and $g \preceq f$.
The equivalence class of 
$f$ is called the $growth\; rate$ of $f$. Simple verification shows that growth functions of 
an algebra with respect 
to different generating sets are equivalent.
\vspace{.2cm}
\\The growth rate is a useful invariant for
finitely generated algebraic structures such as groups, semigroups and algebras.
 The notion of growth function for groups was introduced by Schwarz 
\cite{schwarz55} and independently by Milnor \cite{milnor68}. The description of groups of polynomial
growth was obtained by Gromov in his celebrated work \cite{gromov81}. He proved that every finitely
generated group of polynomial growth contains a nilpotent subgroup of finite index.
\vspace{.2cm}
\par The study of growth of algebras dates back to the papers by Gelfand and Kirillov, \cite{GK661,GK66}.
In this paper
we are mainly interested in finitely presented algebras whose growth functions behave
in intermediate way 
i.e., they grow faster than any polynomial function but slower than any
exponential function. Govorov gave the first examples of finitely generated semigroups
and associative algebras of intermediate growth in \cite{gov72}.
Examples of algebras of intermediate growth can also be found in  
\cite{stewart75,smith76, shearer,ufna80, kril83}. The first examples of finitely generated
groups of intermediate growth were constructed by Grigorchuk \cite{gri83,gri84}.
It is still an open problem
whether 
there exists a finitely presented group of intermediate growth. In contrast, there are examples of
finitely presented algebras of intermediate growth.
The first example is the universal enveloping algebra of a
Lie algebra $W$ with basis 
$\{w_{-1},w_0,w_1,w_2,\dots\}$ and brackets defined by $[w_i,w_j]=(i-j)w_{i+j}$. $W$ is a subalgebra 
of the 
generalized Witt algebra $W_{\mathbb{Z}}$ (see \cite[p.206]{as74} for definitions).
It was proven in \cite{stewart75} that $W$ has a finite presentation with two generators and 
six relations.  It is also a graded algebra with generators of degree $-1$ and $2$. Since $W$ has linear
growth, its universal enveloping algebra is an example of finitely presented associative algebra of 
intermediate growth.
 
 \vspace{.2cm}

 The main goal of this paper is to present examples of finitely presented quadratic algebras
(algebras defined by quadratic relations) of
intermediate growth. The class of quadratic algebras contains a class of finitely presented algebras,
called \emph{Koszul algebras}.
They play an important role in many studies. 
In \cite{PP71}, it is conjectured that the Hilbert series of a Koszul algebra $A$ is
a rational function and
in particular, the growth of $A$ is either polynomial or exponential.
\vspace{.2cm}
 \\In order to construct our first example of 
 a finitely presented quadratic algebra of intermediate growth, we consider 
the Kac-Moody
algebra for the generalized Cartan matrix
$A=\left (\begin{array}{cc}
              2&-2\\
              -2&2
             \end{array} \right)$. This is a
graded Lie algebra of polynomial growth whose generators are of degree $1$. Next, we consider a suitable
subalgebra and its universal enveloping algebra.

\begin{thm}
\label{thm1}
Let  $U$ be 
 the associative algebra with generators $x$, $y$ and relations
 $x^3y-3x^2yx+3xyx^2-yx^3=0$, 
  $y^3x-3y^2xy+3yxy^2-xy^3=0$. Then
  \begin{itemize}
    \item[(i)]  It is the universal enveloping algebra of a subalgebra of 
   the the Kac-Moody
algebra for the generalized Cartan matrix
$A=\left (\begin{array}{cc}
              2&-2\\
              -2&2
             \end{array} \right)$.
  \item[(ii)] $U$ is a graded algebra with generators of degree $1$.
 \item[(iii)] It has intermediate growth of type $e^{\sqrt{n}}$.
\item[(iv)] The Veronese subalgebra $V_4(U)$ of $U$ is a quadratic algebra given by $14$
generators and $96$ quadratic
relations and it has the same growth type with $U$.
  \end{itemize}

\end{thm}
The Kac-Moody algebra for the generalized Cartan matrix
$\di A=\left (\begin{array}{cc}
              2&-2\\
              -2&2
             \end{array} \right)$ is the affine Lie algebra $A_1^{(1)}$. 
             (For the definition of Kac-Moody algebras and classification of affine Lie 
             algebras see \cite{kac}). It has a subalgebra which is 
             isomorphic to the
Lie subalgebra $L$ 
of $sl_{2}(\mathbb{C}[t])$
which consists of all matrices with entries on and under the diagonal divisible 
by t. That is,
$$L=\{a=(a_{ij})_{2\times2}\;|\;\;a_{ij}\in \mathbb{C}[t],\; tr(a)=0\;\text{and for}\;(i,j)\neq(1,2),\; t\;\text{divides}\;a_{ij}\}$$
with the usual Lie bracket $[a,b]=ab-ba$.  It follows from 
 \cite [Theorem 9.11] {kac} that $L$ is finitely presented.
 \vspace{.2cm}
 \\In this paper we will prove this 
 by using the axioms of Lie bracket without mentioning the theory of Kac-Moody algebras.
 In Section \ref{2} we show that $L$ is a finitely presented graded Lie algebra whose generators are all of 
 degree $1$ and $L$ has linear growth.
 In Section \ref{3} we explain the relation between the growth of a Lie algebra and its universal enveloping 
 algebra. In Section \ref{4} we consider
 the Veronese subalgebra of $U$ to obtain a finitely presented quadratic algebra of intermediate 
 growth and in Section \ref{5} we
 complete the proof of Theorem 1. In Section \ref{6} we give another example of finitely presented
 associative 
 algebra $A$ of intermediate growth related to the example of the monoid in \cite{koba}. 
 $A$ has the following presentation:
 $$ A=\langle a,b,c \; \mid \;b^2a=ab^2,\;b^2c=aca,\;acc=0,\;aba=0,\;abc=0,\;cba=0,\;cbc=0  \rangle$$
 We show that $A$ has intermediate growth of type $e^{\sqrt{n}}$ and its Veronese subalgebra $V_3(A)$
 is an example of finitely presented quadratic algebra of intermediate growth.
 In Section \ref{7} , we give an explicit presentation of the Veronese subalgebra $V_4(U)$ 
 of the first construction $U$ as an example of a finitely presented quadratic algebra of intermediate growth.
\vspace{.2cm}

 \section{An example of a finitely presented Lie Algebra of linear growth}
 \label{2}
 
 The following example is a subalgebra of the Kac-Moody Algebra for the generalized Cartan matrix
 $\di A=\left (\begin{array}{cc}
              2&-2\\
              -2&2
             \end{array} \right)$
 \cite{kac}. 
 \smallskip
 
 Consider the subalgebra $L$ of $Sl_{2}(\mathbb{C}[t])$ over $\mathbb{C}$
 (i.e., matrices of trace $0$ with entries in $\mathbb{C}[t])$)
which consists of matrices whose entries on and under the diagonal are divisible by t.
That is, $$L=\{a=(a_{ij})_{2x2}|\;\;a_{ij}\in \mathbb{C}[t],\;
tr(a)=0\;\text{and for}\;(i,j)\neq(1,2),\; t\;\text{divides}\;a_{ij}\}$$
with the usual Lie bracket $[a,b]=ab-ba$. 
\begin{pro}Let $L$ be the Lie algebra described above. Then it has the following properties.

\begin{itemize}
  \item[(i)] $L$ is finitely presented  with generators
  $ x:=
 \left( \begin{array}{cc}
0 & 1  \\
0 & 0 
\end{array} \right)$ and $ y:=\left (\begin{array}{cc}
0&0\\
t&0
\end{array} \right)$ and the defining relations
$[x,[x,[x,y]]]=0$ and 
  $[y,[y,[y,x]]]=0$.
  \item[(ii)] $\di L=\bigoplus_{k\geq 1}L_k$ is graded
  and generated by $L_1$.
  \item[(iii)] L has linear growth.
 \end{itemize}

\end{pro}

\begin{proof}\par
Take $\di x_1:= x=
 \left( \begin{array}{cc}
0 & 1  \\
0 & 0  
 
\end{array} \right) $,
$\di y_1:=y=\left (\begin{array}{cc}
0&0\\
t&0
\end{array} \right)$, and let 
$\di z_1:=\left (\begin{array}{cc}
              t&0\\
              0&-t
             \end{array} \right)
$.
\\In fact, define $\di x_i:=
 \left( \begin{array}{cc}
0 & t^{i-1}  \\
0 & 0  
 
\end{array} \right) $,
$\di y_i:=\left (\begin{array}{cc}
0&0\\
t^i&0
\end{array} \right)$, and let 
$\di z_i:=\left (\begin{array}{cc}
              t^i&0\\
              0&-t^i
             \end{array} \right)
$ for $i\geq 1$.
\\An arbitrary element $w\in L$ is of the form:
$$w=\left( \begin{array}{cc}
            \sum_{i=1}^{n}m_it^i& \sum_{i=1}^{n}k_it^{i-1}\\
             \sum_{i=1}^{n}l_it^i& \sum_{i=1}^{n}-m_it^{i}
           \end{array} \right)=\sum_{i=1}^{n}k_ix_i+ \sum_{i=1}^{n}l_iy_{i}
            + \sum_{i=1}^{n}m_iz_i
$$
So, any element of $L$ can be written as a linear combination of $x_i,\;y_i,\;z_i$ for 
$i\geq 1$ and
 $\di \{x_i,\;y_i,\;z_i\}_{i=1}^{\infty}$ forms a 
linearly independent set
over 
$\mathbb{C}$.
$L$ has the following relations
  \begin{equation}\label{eq:1}[x_i,y_j]=z_{i+j-1},
  \end{equation}
  \begin{equation}\label{eq:2}
 [x_i,z_j]=-2x_{i+j},
 \end{equation}
 \begin{equation}\label{eq:3}
 [y_i,z_j]=2y_{i+j},
 \end{equation}
 \begin{equation}\label{eq:4}
 [x_i,x_j]=0,
 \end{equation}
 \begin{equation}\label{eq:5}
[y_i,y_j]=0,
\end{equation}
\begin{equation}\label{eq:6}
[z_i,z_j]=0.
\end{equation} for $i,j\geq 1$.
 In particular,
 $$x_{i+1}=\di -\frac{1}{2}[x_i,z_1],$$
 $$y_{i+1}=\di \frac{1}{2}[y_i,z_1],$$
 $$z_i=[x_i,y_1].$$
 It follows that $L$ is generated by $x_1$ and $y_1$.
 In order to show that 
 all the relations
 \eqref{eq:1} - \eqref{eq:6} can be derived from the relations 
 $[x_1,[x_1,[x_1,y_1]]]=0$ and 
  $[y_1,[y_1,[y_1,x_1]]]=0$, we apply induction on $i+j=n$.
  If $i+j=2$, 
  the relations \eqref{eq:1} - \eqref{eq:6} hold trivially.
  If $i+j=3$, 

  $$ \di \begin{array}{lll}
       \di  [{x_1},{y_2}]&=&[{x_1},\frac{[{y_1},{z_1}]}{2}]\\
                                  &=&-\frac{1}{2}\left([{z}_1,[{x}_1,{y}_1]]+[{y}_1,[{z}_1,{x}_1]]\right)
                                  \\
                                  &=&[{x}_2,{y}_1]\\
                                  &=&{z}_2
  \end{array}             
$$

\medskip

$$\begin{array}{llll}
    [x_1,z_2]&=&[x_1,[x_2,y_1]]&\\
   &=&-[y_1,[x_1,x_2]]+[x_2,[y_1,x_1]]&\;(\text{since}\; [x_1,x_2]=0)\\
                             &=&[x_2,[x_1,y_1]]&\\
                              &=& [x_2,z_1]&\\
                              &=& -2x_3&
 \end{array}
$$

\medskip

$$\begin{array}{llll}
    [y_1,z_2]&=&[y_1,[x_1,y_2]]&\\
                             &=&-([y_2,[y_1,x_1]]+[x_1,[y_2,y_1]])&\;(\text{since}\; [y_1,y_2]=0)\\
                             &=& [y_2,z_1]&\\
                              &=& 2y_3&
 \end{array}
$$

\medskip
The relations \eqref{eq:4}-\eqref{eq:5} for $n=3$ correspond to relations of $L_0$.
Observe the following three equations for $[z_2,z_1]$,
$$\begin{array}{lll}

[z_2,z_1]&=&[[x_2,y_1],z_1]\\
                         &=&-([[z_1,x_2],y_1]+[[y_1,z_1],x_2])\\
                         &=&[[x_2,z_1],y_1]+[x_2,[y_1,z_1]]\\
                         &=&-2[x_3,y_1]+2[x_2,y_2]\\
                         &=&k

 \end{array}
$$

$$\begin{array}{lll}

[z_2,z_1]&=&[[x_1,y_2],z_1]\\
                         &=&-([[z_1,x_1],y_2]+[[y_2,z_1],x_1])\\
                         &=&[[x_1,z_1],y_2]+[x_1,[y_2,z_1]]\\
                         &=&-2[x_2,y_2]+2[x_1,y_3]\\
                         &=&l

 \end{array}
$$
$$\begin{array}{lll}

[z_2,z_1]&=&[z_2,[x_1,y_1]\\
                         &=&-([y_1,[z_2,x_1]]+[x_1,[y_1,z_2]])\\
                         &=&2[x_3,y_1]-2[x_1,y_3]\\
                         &=&m

 \end{array}
$$
$3.[z_2,z_1]=k+l+m=0$. So, \eqref{eq:1} - \eqref{eq:6} hold for $n=3$.
Now, suppose that
 \eqref{eq:1} - \eqref{eq:6} hold for $i+j\leq n$ for some $n\geq 3$. 
 For $1\leq i\leq n-1$,
$$\di \begin{array}{lll}
   [x_i,y_{j+1}]&=&\frac{1}{2}[x_i,[y_j,z_1]]\\
                                &=&-\frac{1}{2}([z_1,[x_i,y_j]]+[y_j,[z_1,x_i]])\\
                                &=&[x_{i+1},y_{j}]
  \end{array}
$$
\medskip
$$\begin{array}{lll}
   -2x_{n+1}&=&[x_n,z_1]\\
                    &=&-\frac{1}{2}[[x_1,z_{n-1}],z_1]\\
                    &=&\frac{1}{2}([[z_1,x_1],z_{n-1}]+[[z_{n-1},z_1],x_1])\\
                    &=&[x_2,z_{n-1}]
  \end{array}
$$
and,
$$\begin{array}{lll}
   [x_i,z_{j+1}]&=&[x_i,[x_1,y_{j+1}]]\\
                                &=&-([y_{j+1},[x_i,x_1]]+[x_1,[y_{j+1},x_i]])\\
                                &=&[x_1,z_{i+j}]
  \end{array}
$$
 Similarly, it can be shown that 
 
 $$ 2y_{n+1}=[y_i,z_{j+1}]$$
 for any $i,j\geq 1$ such that $i+j=n$.
So \eqref{eq:1} - \eqref{eq:3} hold for $i+j=n+1.$ 

$$\begin{array}{lll}
[x_1,x_n]&=&-\frac{1}{2}[x_1,[x_i,z_j]]\\
                         &=&\frac{1}{2}([z_j,[x_1,x_i]]+[x_i,[z_j,x_1]])\\
                         &=&-\frac{1}{2}[x_i,[x_1,z_j]]\\
                         &=&[x_i,x_j]\\
                        
                      \end{array}
$$
This equality implies
$[x_i,x_j]=[x_j,x_i]$. Similarly, one checks that 
$[y_i,y_j]=[y_j,y_i]$. Hence, \eqref{eq:4} - \eqref{eq:5} hold for $i+j=n+1.$

\medskip
Finally, we need check that \eqref{eq:6}  holds for $i+j=n+1$. 

$$\begin{array}{lllll}
   [z_1,z_n]&=&[z_1,[x_n,y_1]]&=&2[x_{n+1},y_1]-2[x_n,y_2]\\
                            &=&[z_1,[x_{n-1},y_2]]&=&2[x_{n},y_2]-2[x_{n-1},y_3]\\
                            & &\vdots                                    &&\\
                             &=&[z_1,[x_{1},y_n]]&=&2[x_{2},y_{n}]-2[x_{1},y_{n+1}]
                              
  \end{array}
$$
implies that $n.[z_1,z_n]=2[x_{n+1},y_1]-2[x_{1},y_{n+1}]$ and,
$$\begin{array}{lllll}
   2[x_{1},y_{n+1}]&=&[x_1,[y_1,z_n]]&=
   &-[z_n,[x_1,y_1]]-[y_1,[z_n,x_1]]\\
   &&&=&[z_1,z_n]+2[x_{n+1},y_1]
  \end{array}
$$
So $[z_1,z_n]=0$.
Now, consider $[z_i,z_j]$ for $i\in\{1,\dots,n-1\}$,
$$\begin{array}{lllll}
   [z_i,z_j]&=&[z_i,[x_j,y_1]]&=&-([y_1,[z_i,x_j]]+[x_j,[y_1,z_i]])\\
                            & &          &= &2[x_{i+j},y_1]-2[x_j,y_{i+1}] 
  \end{array}
  $$
and,
$$\begin{array}{lllll}
   [x_j,y_{i+1}]&=&\frac{1}{2}[x_j,[y_i,z_1]&=&-\frac{1}{2}([z_1,[x_j,y_i]]+
   [y_i,[z_1,x_j]])\\
                                & &                  &=&-\frac{1}{2}([z_1,z_n]+[y_i,2x_{j+1}])\\
                                & &                                                 &=&[x_{j+1},y_i]
  \end{array}
$$
By applying this $i$ times we get $[x_j,y_{i+1}]=[x_n,y_1]$ , so that
$$[z_i,z_j]=0\;\text{for}\;i+j=n+1$$ i.e., \eqref{eq:6}  holds for $i+j=n+1$.
By \eqref{eq:1} - \eqref{eq:3}, the set
$\di \{x_i,y_i,z_i\}_{i=1}^{\infty}$
forms a basis for $L$ as a vector space. It can be observed that
$L=\di \bigoplus_{k\geq 1}L_k$ where 
$L_{2k-1}=\langle x_k\rangle \oplus \langle y_k \rangle$
and $L_k=\langle z_k\rangle$ for $k\geq 1$. Since
$$[L_{2k-1},L_{2m-1}]\subseteq L_{2(k+m-1)},$$ 
$$[L_{2k},L_{2m}]=0,$$
$$[L_{2k-1},L_{2m}]\subseteq L_{2(k+m)-1},$$
$L$ admits an $\mathbb{N}$-gradation given by the sum of occurrences of
$x$ and $y$ in each commutator i.e.,
$ L=\bigoplus_{k\geq 1}L_k$ is a graded Lie algebra generated
by two elements of degree $1$ 
$(deg(a)= min\{n| a\in \bigoplus_{k= 1}^{n}L_k)\})$ and $L$ has
linear growth ($dim\;L_i\in \{1,2\}$ for $i\geq 1$ ).

\end{proof}

\begin{rem}
We notice that $L$ also admits a $\mathbb{Z}$-gradation. It is a 3-graded Lie algebra 
(in the sense of \cite{do}) over $\mathbb{C}$ 
generated by elements  $x$ of degree $1$ and $y$ of degree $-1$ .

\end{rem}

\medskip

\section{The relation between the growth of a Lie algebra and its universal enveloping algebra} 
\label{3}
 
     Let $L$ be any Lie algebra over  a field $k$ and $U(L)$ be its 
     universal enveloping algebra. For an ordered basis 
 $u_1,u_2,\dots$
 of $L$, monomials $u_{i_1}\dots u_{i_r}$ with $i_1\leq i_2 \leq \dots \leq i_r$
 form a basis for $U(L)$ (Poincar\'e-Birkhoff-Witt Theorem (\cite {diamond})).
If $L=\bigoplus {L_n}$ is a graded Lie algebra such that all the components are finite dimensional,
then
\begin{equation}\label{eq:7}
\displaystyle \sum_{n=0}^{\infty} b_n t^n =\prod_{n=1}^{\infty}(1-t^n)^{-a_n}
\end{equation}
where $a_n:=dim(L_n)$ and $b_n$:=number of monomials of length $n$ in $U(L)$ (\cite{smith76}).
 The proof of the following proposition can be found in
 various papers (\cite{ber83}, \cite{pet93}, \cite{gribar00}).
\begin{pro}
\label{pro2}
 If $a_n$ and $b_n$ are related by \eqref{eq:7} and
 $a_n\sim n^d$, then 
 $b_n \sim e^{n^{\frac{d+1}{d+2}}}$.
  
 \end{pro}

\begin{cor}
If a Lie algebra $L$ grows polynomially then its universal enveloping algebra $U(L)$
 has intermediate growth. 
 In particular, if $L$ has linear growth, then $U(L)$ has growth of type $e^{\sqrt{n}}$.
 \end{cor}
 
 \section{Veronese subalgebra of an associative graded algebra}
 \label{4}
  Let $A=k\langle x_1,\dots,x_m\rangle$ be a free associative algebra over a field $k$
  with generating set $\{x_1,\dots,x_m\} $.
Each element $u$ of $A$ can be written uniquely
   as 
   $$u=u_0+u_1+\dots +u_l ,$$
   where $A_0=k$, $u_i \in A_i$  and $A_i$ is the vector space over $k$ spanned by $m^i$ monomials of length $i$.
    Let $R =\{ f_1, f_2, \dots, f_s\}$ be a finite set of non-zero  homogeneous polynomials
   and $I$ be the ideal generated by
   $R$. 
  Since $I$ is generated by homogeneous polynomials, the factor algebra
  $\tilde{A}=A/I$ is graded:
  $$\tilde{A}=\tilde{A_0}\oplus\tilde{A_1}\oplus \dots \oplus\tilde{ A_n}\oplus \dots $$ 
  where $\tilde{A_i}=(A_i+I)/I\cong A_i/(A_i \cap I)$. For $d\geq 1$, a \emph{Veronese subalgebra}
  of $\tilde{A}$ is defined as
  $$V_d(\tilde{A}):= k \oplus \tilde{A}_d \oplus \tilde{A}_{2d}\oplus \dots$$
  It is straightforward to see that,
  $$growth\;of\; \tilde{A}\sim growth\;of\;V_d(\tilde{A})$$
  \begin{pro} \cite{bac85}
  For sufficiently large $d$, $V_d(\tilde{A})$ is quadratic.
  \end{pro}
  \begin{proof}
   Let $d_1,\dots, d_s$ be the degrees of  $f_1, f_2, \dots, f_s$ respectively and
  $d\geq max\{d_i,\;1\leq i \leq s\}$. For any two words $v',\; v''$
  such that 
  $$deg(v')+d_i+deg(v'')=d$$
  consider the element $v'f_i v''\in A_d$, and for any two words $w',\; w''$ such that
  $$deg(w')+d_i+deg(w'')=2d$$
  consider the element $w'f_i w''\in A_{2d}$. Let $R^{*}=\{v'f_i v'', w'f_i w''\}$
  for $i \in \{ 1,\dots, s\}$ and $a$ be a homogeneous element from  $A^{(n)}\cap I$.
  Say $a=\sum \alpha v f_i w$, 
 where $\alpha \in k$, $v$ and $w$ are words. If we choose a summand
 and represent $v=v_1 v_2$, $deg(v_1)$ is a multiple of $d$, $0\leq deg(v_2) < d $.
 Similarly, $w=w_2w_1$, $deg(w_1)$ is a multiple of $d$, 
 $0 \leq deg(w_2) < d$. Then we will get $deg(v_2 f_i w_2)=d$ or $2d$. 
 Hence $v_2f_i w_2 \in R^{*}$. It shows that $V_d(A)\cap I$ is an ideal generated by the elements of $R^*$
 and an element $v'f_i v''$ is a linear combination of free generators of $A^{(n)}$ whereas
 $w'f_i w''$ is a quadratic element in these generators. So $V_d(\tilde{A})=V_d(A)/(V_d(A)\cap I)$ is a
 quadratic algebra.
   
  \end{proof}

 \section{Proof of Theorem 1}
 \label{5}
 Let $\di L=\langle x_1,\dots,x_m \mid f_1=0,\dots, f_r=0\rangle$ where each of $f_i$ 
 is a linear combination of the commutators (elements of the form $[x_{i_1},\dots,x_{i_k} ]$
 with an arbitrary distribution of parentheses inside). Then the 
 universal enveloping algebra $U(L)$
 of $L$ is an associative algebra with the identical set of generators and relations,
 where the commutators
 are thought of as in the ordinary associative sense: $[x,y]=xy-yx$
 \cite[Proposition 2, p.14]{Bou}.
The universal enveloping algebra $U(L)$ of $L=
\langle x_1, y_1\;\mid \;[x_1,[x_1,[x_1,y_1]]]=0,\;
  [y_1,[y_1,[y_1,x_1]]]=0\rangle$ has the following presentation:
$$U(L)=\langle x_1, y_1\; \mid
\; x_1^3y_1-3x_1^2y_1x_1+3x_1y_1x_1^2-y_1x_1^3=0,
\;y_1^3x_1-3y_1^2x_1y_1+3y_1x_1y_1^2-x_1y_1^3=0\rangle.$$
So, the associative algebra $U$ in Theorem \ref{thm1}
 is the universal enveloping algebra $U(L)$ of
 $L$. By Proposition \ref{pro2}, since $L$ has linear growth, the growth rate of $U(L)$
 is intermediate of type
 $e^{\sqrt{n}}$ .
 In order to obtain 
 a quadratic algebra of intermediate growth we consider a Veronese subalgebra of 
 $V_4(U)$ as explained
 in the previous section and
 conclude that for a given finitely presented graded algebra with all generators of degree $1$, 
 one can construct a finitely presented
graded algebra with all relations of degree $2$. $V_4(U)$ is an example of
a finitely presented graded algebra with intermediate growth.
It has $14$ generators and $96$ relations. In the next section we compute all these relations.

\begin{section}{A construction based on Kobayashi's example}
\label{6}
 
In this section we construct another example of a finitely presented associative algebra with quadratic relations whose 
growth function is intermediate. For this, we consider the following example of  a monoid with $0$
that appears in the
paper of Kobayashi \cite{koba}. 
$$M=\langle a,b,c \; \mid \;ba=ab,\;bc=aca,\;acc=0  \rangle$$
where $w(a)=w(c)=1,\; w(b)=2$, $w$ is a positive weight function on $M$. Kobayashi
shows that $M$ is a finitely presented monoid with solvable word problem which cannot be presented by a regular
complete system. In order to prove that it cannot be presented by a regular complete system, he proves that 
$M$ has intermediate growth. Now, we consider the semigroup algebra $k[M]$ over a field $k$. $k[M]$
has the same presentation and growth function with $M$. So $k[M]$ is an example of finitely presented associative graded 
algebra of intermediate growth. But the generators of $k[M]$ have degrees $deg(a)=deg(c)=1$ and 
$deg(b)=2$. To construct a quadratic algebra with these properties, we need to consider an algebra whose
generators are all of degree $1$. Thus we consider the following monoid:
$$ \tilde{M}=\langle a,b,c \; \mid \;b^2a=ab^2,\;b^2c=aca,\;acc=0,\;aba=0,\;abc=0,\;cba=0,\;cbc=0  \rangle$$ where
$w(a)=w(b)=w(c)=1$.
\vspace{.3cm}
\\Now, we have the monoid algebra $A:=k[\tilde{M}]$ over a field $k$:

$$ A=\langle a,b,c \; \mid \;b^2a=ab^2,\;b^2c=aca,\;acc=0,\;aba=0,\;abc=0,\;cba=0,\;cbc=0  \rangle$$ where
$deg(a)=deg(b)=deg(c)=1$
To show that $A$ has intermediate growth, we first find a complete rewriting system for $A$.
Let $\prec$ be the shortlex order on $\langle X \rangle$ based on the order $a\prec b\prec c$
i.e.,
$$\di w_1\prec w_2 \; \text{implies}\; |w_1|<|w_2|\;\; \text{or}\;\; |w_1|=|w_2|
\;\&\; w_1 \prec _{lex} w_2.$$
 Then $A$ has the rewriting system $R $ consisting of the following relations
 $$\begin{array}{ccc}
    b^2a &\rightarrow &ab^2\\
  b^2c &\rightarrow& aca\\
  acc&\rightarrow &0\\
 aba &\rightarrow & 0\\
abc &\rightarrow & 0\\
  cba & \rightarrow & 0\\
  cbc & \rightarrow & 0
   \end{array}
$$
It is easily seen that $R$ is Noetherian. By applying the \emph{Knuth-Bendix algorithm},
we obtain the following complete rewriting system 
 $R_{\infty}$ equivalent to $R$:
 $$\di R_{\infty}=\{b^2a\rightarrow ab^2,\;b^2c\rightarrow aca,\; aba \rightarrow  0,\;
abc \rightarrow  0,\;
  cba  \rightarrow  0,\;
  cbc  \rightarrow  0\}\cup 
 \bigcup_{n=1}^{\infty}\{a^nca^{n-1}c\rightarrow 0\}$$
 
 A monomial (word) $m$ is called \emph{irreducible}  with respect to the rewriting system $R$
 if all the rewriting rules 
 act trivially on $m$. The set of all irreducible words with respect to $R$ is denoted by $Irr(R)$. 
 Since $R_{\infty}$ is a complete rewriting system, $Irr(R_{\infty})$ is the set of words which
 do not contain $u$ as a subword for any  $u\rightarrow v \in R_{\infty}$.
 By \emph{Bergman's Diamond Lemma}
 \cite{diamond}, $Irr(R_{\infty})$, forms a basis for $A$. Words in $Irr(R_{\infty})$ are of the following 
 form
 $$b^sa^{m_1}ca^{m_2}c\dots a^{m_r}ca^{l}b^{k}$$
 where $s\in\{0,1\}$, $l,k\in \mathbb{N}\cup \{0\}$ and $0\leq m_1 \leq m_2 \leq \dots \leq m_r$
 , $m_i \in \mathbb{N}\cup \{0\} $ for $i\in \{1,\dots r\}$.
 So, the number of words in $Irr(R^{\infty})$ of length $n$ is equal to
$$\di \sum_{j=0}^{n}(2j+1)\cdot\lvert \{(m_1,\dots,m_r) \mid 
0\leq m_1\leq\dots\leq m_r,\;m_1+\dots+m_r=n-j-r\}\rvert$$
$$\di =\sum_{j=0}^{n}(2j+1)\cdot p(n-j)$$
where $p(n)$ is the number of partitions of $n$. Hence
$$\di \gamma_A(n)\sim p(n)\sim e^{\sqrt{n}}.$$
 $A$ is an example of finitely presented graded algebra 
with generators of degree $1$ and intermediate growth function and its 
 \emph {Veronese subalgebra} $V_3(A)$ can be presented by finitely many quadratic relations 
 (to be precise with $21$ generators
 and $280$ relations).
 
 \end{section}

\begin{section}{Appendix: Presentation of the Veronese subalgebra $V_4(U)$ of $U$}
\label{7}
As we noted in the previous section, $U(L)$ is an associative algebra with generators $x,y$
and the set of relations $R=\{x^3y-3x^2yx+3xyx^2-yx^3=0,
\;y^3x-3y^2xy+3yxy^2-xy^3=0\}$. Since $R$ is a set of two homogeneous polynomials,
$U$ is a graded algebra.
Let $V_4(U)$ be the Veronese subalgebra of $U$.
It was proven in Section \ref{4} that $V_4(U)$ is a graded algebra
generated by the set 
$S$ of monomials of length $4$ over $\{x,y\}$ and the set of relations 
$R^{*}=\{f_i=0,vf_iw=0\}$ where $v,w$ are monomials such that $l(v)+l(w)=4$ and, 
$f_1=x^3y-3x^2yx+3xyx^2-yx^3,
\;f_2=y^3x-3y^2xy+3yxy^2-xy^3.$
Basically, $R^{*}$ is the set of homogeneous polynomials of degree $4$ or $8$ generated by 
$R=\{f_1=0,f_2=0\}$ in $k[x,y]$. Since there are $48$ different pairs $(v,w)$ of monomials,
$R^*$ consists of $2$ homogeneous polynomials of degree $4$:
$$(i)\;yx^3=x^3y-3x^2yx+3xyx^2,\;\;(ii)\;y^3x=xy^3-3yxy^2+3y^2xy$$
and $96$
homogeneous polynomials of degree $8$:
 \begin{small}
$$
\begin{array}{ll}

(1)\;xyx^2x^4=x^4yx^3-3x^3yx^4+3x^2yxx^4,& (49)\;x^2yxx^4=x^4xyx^2-3x^4yx^3+3x^3yx^4,\\
(2)\;x^3yx^4=x^4x^2yx-3x^4xyx^2+3x^4yx^3,& (50)\;xy^3x^4=x^2y^2yx^3-3xyxyyx^3+3xy^2xyx^3,\\
(3)\;x^2y^2yx^3=x^3yy^2x^2-3x^2yxy^2x^2+3x^2y^2xyx^2, &(51)\; x^3yy^2x^2=x^4y^3x-3x^3yxy^2x+3x^3yyxyx,\\
(4)\;xyx^2x^3y=x^4yx^2y-3x^3yx^3y+3x^2yxx^3y,&(52)\;x^2yxx^3y=x^4xyxy-3x^4yx^2y+3x^3yx^3y,\\
(5)\;x^3yx^3y=x^4x^2y^2-3x^4xyxy+3x^4yx^2y,& (53)\;xy^3x^3y=x^2y^2yx^2y-3xyxyyx^2y+3xy^2xyx^2y,\\
(6)\;x^2y^2yx^2y=x^3yy^2xy-3x^2yxy^2xy+3x^2y^2xyxy,&(54)\; x^3yy^2xy=x^4y^4-3x^3yxy^3+3x^3yyxy^2,\\
(7)\;xyx^2x^2yx=x^4yxyx-3x^3yx^2yx+3x^2yxx^2yx,&(55)\; x^2yxx^2yx=x^4xy^2x-3x^4yxyx+3x^3yx^2yx,\\
(8)\;x^2y^2x^4=x^2yxx^2yx-3x^2yxxyx^2+3x^2yxyx^3,& (56)\;xy^3x^2yx=x^2y^2yxyx-3xyxyyxyx+3xy^2xyxyx,\\
(9)\;x^2y^2yxyx=x^3yy^3x-3x^2yxy^3x+3x^2y^2xy^2x,& (57)\;x^2y^2y^2x^2=x^2yxy^3x-3x^2y^2xy^2x+3x^2y^2yxyx,\\
(10)\;xyx^2x^2y^2=x^4yxy^2-3x^3yx^2y^2+3x^2yxx^2y^2,&(58)\; x^2yxx^2y^2=x^4xy^3-3x^4yxy^2+3x^3yx^2y^2,\\
(11)\;x^2y^2x^3y=x^2yxx^2y^2-3x^2yxxyxy+3x^2yxyx^2y,& (59)\;xy^3x^2y^2=x^2y^2yxy^2-3xyxyyxy^2+3xy^2xyxy^2,\\
(12)\;x^2y^2yxy^2=x^3yy^4-3x^2yxy^4+3x^2y^2xy^3,&(60)\; x^2y^2y^2xy=x^2yxy^4-3x^2y^2xy^3+3x^2y^2yxy^2,\\
(13)\;xyx^2xyx^2=x^4y^2x^2-3x^3yxyx^2+3x^2yxxyx^2,&(61)\; xy^2xx^4=xyx^2xyx^2-3xyx^2yx^3+3xyxyx^4,\\
(14)\;xyxyx^4=xyx^2x^2yx-3xyx^2xyx^2+3xyx^2yx^3,&(62)\; xy^3xyx^2=x^2y^2y^2x^2-3xyxyy^2x^2+3xy^2xy^2x^2,\\
(15)\;xy^3yx^3=xyxyy^2x^2-3xy^2xy^2x^2+3xy^3xyx^2,&(63)\; xyxyy^2x^2=xyx^2y^3x-3xyxyxy^2x+3xyxyyxyx,\\
(16)\;xyx^2xy^2x=x^4y^3x-3x^3yxy^2x+3x^2yxxy^2x,&(64)\; xy^2xx^2yx=xyx^2xy^2x-3xyx^2yxyx+3xyxyx^2yx,\\
(17)\;xy^3x^4=xy^2xx^2yx-3xy^2xxyx^2+3xy^2xyx^3,&(65)\; xy^3xy^2x=x^2y^2y^3x-3xyxyy^3x+3xy^2xy^3x,\\
(18)\;xy^3yxyx=xyxyy^3x-3xy^2xy^3x+3xy^3xy^2x,&(66)\; xy^3y^2x^2=xy^2xy^3x-3xy^3xy^2x+3xy^3yxyx,\\
(19)\;xyx^2xyxy=x^4y^2xy-3x^3yxyxy+3x^2yxxyxy,&(67)\; xy^2xx^3y=xyx^2xyxy-3xyx^2yx^2y+3xyxyx^3y,\\
(20)\;xyxyx^3y=xyx^2x^2y^2-3xyx^2xyxy+3xyx^2yx^2y,&(68)\; xy^3xyxy=x^2y^2y^2xy-3xyxyy^2xy+3xy^2xy^2xy,\\
(21)\;xy^3yx^2y=xyxyy^2xy-3xy^2xy^2xy+3xy^3xyxy,&(69)\; xyxyy^2xy=xyx^2y^4-3xyxyxy^3+3xyxyyxy^2,\\
(22)\;xyx^2xy^3=x^4y^4-3x^3yxy^3+3x^2yxxy^3,&(70)\;xy^2xx^2y^2=xyx^2xy^3-3xyx^2yxy^2+3xyxyx^2y^2,\\
(23)\;xy^3x^3y=xy^2xx^2y^2-3xy^2xxyxy+3xy^2xyx^2y,&(71)\; xy^3xy^3=x^2y^2y^4-3xyxyy^4+3xy^2xy^4,\\
(24)\;xy^3yxy^2=xyxyy^4-3xy^2xy^4+3xy^3xy^3,&(72)\; xy^3y^2xy=xy^2xy^4-3xy^3xy^3+3xy^3yxy^2,\\
(25)\;y^2x^2x^4=yx^3yx^3-3yx^2yx^4+3yxyxx^4,&(73)\; yxyxx^4=yx^3xyx^2-3yx^3yx^3+3yx^2yx^4,\\
(26)\;yx^2yx^4=yx^3x^2yx-3yx^3xyx^2+3yx^3yx^3,&(74)\; y^4x^4=yxy^2yx^3-3y^2xyyx^3+3y^3xyx^3,\\
(27)\;yxy^2yx^3=yx^2yy^2x^2-3yxyxy^2x^2+3yxy^2xyx^2,&(75)\; yx^2yy^2x^2=yx^3y^3x-3yx^2yxy^2x+3yx^2yyxyx,\\
(28)\;x^2y^2x^2yx=yx^3yxyx-3yx^2yx^2yx+3yxyxy^2xy,&(76)\; yxyxx^2yx=yx^3xy^2x-3yx^3yxyx+3yx^2yx^2yx,\\
(29)\;yxy^2x^4=yxyxx^2yx-3yxyxxyx^2+3yxyxyx^3,&(77)\; y^4x^2yx=yxy^2yxy x-3y^2xyyxyx+3y^3xyxyx,\\
(30)\;yxy^2yxyx=yx^2yy^3x-3yxyxy^3x+3yxy^2xy^2x,&(78)\; yxy^2y^2x^2=yxyxy^3x-3yxy^2xy^2x+3yxy^2yxyx,\\
(31)\;y^2x^2x^2y^2=yx^3yxy^2-3yx^2yx^2y^2+3yxyxx^2y^2,&(79)\;yxyxx^2y^2=yx^3xy^3-3yx^3yxy^2+3yx^2yx^2y^2,\\
(32)\;yxy^2x^3y=yxyxx^2y^2-3yxyxxyxy+3yxyxyx^2y,&(80)\; y^4x^2y^2=yxy^2yxy^2-3y^2xyyxy^2+3y^3xyxy^2,\\
(33)\;yxy^2yxy^2=yx^2yy^4-3yxyxy^4+3yxy^2xy^3,&(81)\; yxy^2y^2xy=yxyxy^4-3yxy^2xy^3+3yxy^2yxy^2,\\
(34)\;y^2x^2x^3y=yx^3yx^2y-3yx^2yx^3y+3yxyxx^3y,&(82)\; yxyxx^3y=yx^3xyxy-3yx^3yx^2y+3yx^2yx^3y,\\
(35)\;yx^2yx^3y=yx^3x^2y^2-3yx^3xyxy+3yx^3yx^2y,&(83)\; y^4x^3y=yxy^2yx^2y-3y^2xyyx^2y+3y^3xyx^2y,\\
(36)\;yxy^2yx^2y=yx^2yy^2xy-3yxyxy^2xy+3yxy^2xyxy,&(84)\; yx^2yy^2xy=yx^3y^4-3yx^2yxy^3+3yx^2yyxy^2,\\
(37)\;y^2x^2xyx^2=yx^3y^2x^2-3yx^2yxyx^2+3yxyxxyx^2,&(85)\; y^3xx^4=y^2x^2xyx^2-3y^2x^2yx^3-3y^2xyx^4,\\
(38)\;y^2xyx^4=y^2x^2x^2yx-3y^2x^2xyx^2+3y^2x^2yx^3,&(86)\; y^4xyx^2=yxy^2y^2x^2-3y^2xyy^2x^2+3y^3xy^2x^2,\\
(39)\;y^4yx^3=y^2xyy^2x^2-3y^3xy^2x^2+3y^4xyx^2,&(87)\; y^2xyy^2x^2=y^2x^2y^3x-3y^2xyxy^2x+3y^2xyyxyx,\\
(40)\;y^2x^2xyxy=yx^3y^2xy-3yx^2yxyxy+3yxyxxyxy,&(88)\; y^3xx^3y=y^2x^2xyxy-3y^2x^2yx^2y+3y^2xyx^3y,\\
(41)\;y^2xyx^3y=y^2x^2x^2y^2-3y^2x^2xyxy+3y^2x^2yx^2y,&(89)\; y^4xyxy=yxy^2y^2xy-3y^2xyy^2xy+3y^3xy^2xy,\\
(42)\;y^4yx^2y=y^2xyy^2xy-3y^3xy^2xy+3y^4xyxy,&(90)\; y^2xyy^2xy=y^2x^2y^4-3y^2xyxy^3+3y^2xyyxy^2,\\
(43)\;y^2x^2xy^2x=yx^3y^3x-3yx^2yxy^2x+3yxyxxy^2x,&(91)\; y^3xx^2yx=y^2x^2yx^2x-3y^2x^2yxyx+3y^2xyx^2yx,\\
(44)\;y^4x^4=y^3xx^2yx-3y^3xxyx^2+3y^3xyx^3,&(92)\; y^4xy^2x=yxy^2y^3x-3y^2xyy^3x+3y^3xy^3x,\\
(45)\;y^4yxyx=y^2xyy^3x-3y^3xy^3x+3y^4xy^2x,&(93)\; y^4y^2x^2=y^3xy^3x-3y^4xy^2x+3y^4yxyx,\\
(46)\;y^2x^2xy^3=yx^3y^4-3yx^2yxy^3+3yxyxxy^3,&(94)\; y^3xx^2y^2=y^2x^2xy^3-3y^2x^2yxy^2+3y^2xyx^2y^2,\\
(47)\;y^4x^3y=y^3xx^2y^2-3y^3xxyxy+3y^3xyx^2y,&(95)\; y^4xy^3=yxy^2y^4-3y^2xyy^4-3y^2xyy^4+3y^3xy^4,\\
(48)\;y^4yxy^2=y^2xyy^4-3y^3xy^4+3y^4xy^3,&(96)\; y^4y^2xy=y^3xy^4-3y^4xy^3+3y^4yxy^2.

\end{array}
$$
\end{small}

We can rename the generators as follows:
$$y^4=Y_1,\;y^3x=Y_2,\;y^2xy=Y_3,\;y^2x^2=Y_4,\;yxy^2=Y_5,\;yxyx=Y_6,yx^2y=Y_7,\;yx^3=Y_8,$$
$$xy^3=X_1,\;xy^2x=X_2,\;xyxy=X_3,\;xyx^2=X_4,\;x^2y^2=X_5,\;x^2yx=X_6,\;x^3y=X_7,\;x^4=X_8.$$
So the relations will be
$$(i)\;Y_8=X_7-3X_6+3X_4,\;\; (ii)\; Y_2=X_1-3Y_5+3Y_3$$
\begin{small}
$$
\begin{array}{ll}
(1)\;X_4X_8=X_8Y_8-3X_7X_8+3X_6X_8, & (49)\;X_6X_8=X_8X_4-3X_8Y_8+3X_7X_8,\\
(2)\;X_7X_8=X_8X_6-3X_8X_4+3X_8Y_8, & (50)\;X_1X_8=X_5Y_8-3X_3Y_8+3X_2Y_8,\\
(3)\;X_5Y_8=X_7Y_4-3X_6Y_4+3X_5X_4, & (51)\;X_7Y_4=X_8Y_2-3X_7X_2+3X_7Y_6,\\
(4)\;X_4X_7=X_8Y_7-3X_7X_7+3X_6X_7, & (52)\;X_6X_7=X_8X_3-3X_8Y_7+3X_7X_7,\\
(5)\;X_7X_7=X_8X_5-3X_8X_3+3X_8Y_7, & (53)\;X_1X_7=X_5Y_7-3X_3Y_7+3X_2Y_7,\\
(6)\;X_5Y_7=X_7Y_3-3X_6Y_3+3X_5X_3, & (54)\;X_7Y_3=X_8Y_1-3X_7X_1+3X_7Y_5,\\
(7)\;X_4X_6=X_8Y_6-3X_7X_6+3X_6X_6, & (55)\;X_6X_6=X_8X_2-3X_8Y_6+3X_7X_6,\\
(8)\;X_5X_8=X_6X_6-3X_6X_4+3X_6Y_8, & (56)\;X_1X_6=X_5Y_6-3X_3Y_6+3X_2Y_6,\\
(9)\;X_5Y_6=X_7Y_2-3X_6Y_2+3X_5X_2, & (57)\;X_5Y_4=X_6Y_2-3X_5X_2+3X_5Y_6,\\
(10)\;X_4X_5=X_8Y_5-3X_7X_5+3X_6X_5, & (58)\;,X_6X_5=X_8X_1-3X_8Y_5+3X_7X_5\\
(11)\;X_5X_7=X_6X_5-3X_6X_3+3X_6Y_7, & (59)\;X_1X_5=X_5Y_5-3X_3Y_5+3X_2Y_5,\\
(12)\;X_5Y_5=X_7Y_1-3X_6Y_1+3X_5X_1, & (60)\;X_5Y_3=X_6Y_1-3X_5X_1+3X_5Y_5,\\
(13)\;X_4X_4=X_8Y_4-3X_7X_4+3X_6X_4, & (61)\;X_2X_8=X_4X_4-3X_4Y_8+3X_3X_8,\\
(14)\;X_3X_8=X_4X_6-3X_4X_4+3X_4Y_8, & (62)\;X_1X_4=X_5Y_4-3X_3Y_4+3X_2Y_4,\\
(15)\;X_1Y_8=X_3Y_4-3X_2Y_4+3X_1X_4, & (63)\;X_3Y_4=X_4Y_2-3X_3X_2+3X_3Y_6,\\
(16)\;X_4X_2=X_8Y_2-3X_7X_2+3X_6X_2, & (64)\;X_2X_6=X_4X_2-3X_4Y_6+3X_3X_6,\\
(17)\;X_1X_8=X_2X_6-3X_2X_4+3X_2Y_8, & (65)\;X_1X_2=X_5Y_2-3X_3Y_2+3X_2Y_2,\\
(18)\;X_1Y_6=X_3Y_2-3X_2Y_2+3X_1X_2, & (66)\;X_1Y_4=X_2Y_2-3X_1X_2+3X_1Y_6,\\
 (19)\;X_4X_3=X_8Y_3-3X_7X_3+3X_6X_3, & (67)\;X_2X_7=X_4X_3-3X_4Y_7+3X_3X_7,\\
 (20)\;X_3X_7=X_4X_5-3X_4X_3+3X_4Y_7, & (68)\;X_1X_3=X_5Y_3-3X_3Y_3+3X_2Y_2,\\
 (21)\;X_1Y_7=X_3Y_3-3X_2Y_3+3X_1X_3, & (69)\;X_3Y_3=X_4Y_1-3X_3X_1+3X_3Y_5,\\
 (22)\;X_4X_1=X_8Y_1-3X_7X_1+3X_6X_1, & (70)\;X_2X_5=X_4X_1-3X_4Y_5+3X_3X_5,\\
 (23)\;X_1X_7=X_2X_5-3X_2X_3+3X_2Y_7, & (71)\;X_1X_1=X_5Y_1-3X_3Y_1+3X_2Y_1,\\
 (24)\;X_1Y_5=X_3Y_1-3X_2Y_1+3X_1X_1, & (72)\;X_1Y_3=X_2Y_1-3X_1X_1+3X_1Y_5,\\
 (25)\;Y_4X_8=Y_8Y_8-3Y_7X_8+3Y_6X_8, & (73)\;Y_6X_8=Y_8X_4-3Y_8Y_8+3Y_7X_8,\\
(26)\;Y_7X_8=Y_8X_6-3Y_8X_4+3Y_8Y_8, & (74)\;Y_1X_8=Y_5Y_8-3Y_3Y_8+3Y_2Y_8,\\
(27)\;Y_5Y_8=Y_7Y_4-3Y_6Y_4+3Y_5X_4, & (75)\;Y_7Y_4=Y_8Y_2-3Y_7X_2+3Y_7Y_6,\\
(28)\;Y_4X_6=Y_8Y_6-3Y_7X_6+3Y_6X_6, & (76)\;Y_6X_6=Y_8X_2-3Y_8Y_6+3Y_7X_6,\\
(29)\;Y_5X_8=Y_6X_6-3Y_6X_4+3Y_6Y_8, & (77)\;Y_1X_6=Y_5Y_6-3Y_3Y_6+3Y_2Y_6,\\
(30)\;Y_5Y_6=Y_7Y_2-3Y_6Y_2+3Y_5X_2, & (78)\;Y_5Y_4=Y_6Y_2-3Y_5X_2+3Y_5Y_6,\\
(31)\;Y_4X_5=Y_8Y_5-3Y_7X_5+3Y_6X_5, & (79)\;Y_6X_5=Y_8X_1-3Y_8Y_5+3Y_7X_5\\
(32)\;Y_5X_7=Y_6X_5-3Y_6X_3+3Y_6Y_7, & (80)\;Y_1X_5=Y_5Y_5-3Y_3Y_5+3Y_2Y_5,\\
(33)\;Y_5Y_5=Y_7Y_1-3Y_6Y_1+3Y_5X_1, & (81)\;Y_5Y_3=Y_6Y_1-3Y_5X_1+3Y_5Y_5,\\
(34)\;Y_4X_7=Y_8Y_7-3Y_7X_7+3Y_6X_7, & (82)\;Y_6X_7=Y_8X_3-3Y_8Y_7+3Y_7X_7,\\
(35)\;Y_7X_7=Y_8X_5-3Y_8X_3+3Y_8Y_7, & (83)\;Y_1X_7=Y_5Y_7-3Y_3Y_7+3Y_2Y_7,\\
(36)\;Y_5Y_7=Y_7Y_3-3Y_6Y_3+3Y_5X_3, & (84)\;Y_7Y_3=Y_8Y_1-3Y_7X_1+3Y_7Y_5,\\
(37)\;Y_4X_4=Y_8Y_4-3Y_7X_4+3Y_6X_4, & (85)\;Y_2X_8=Y_4X_4-3Y_4Y_8+3Y_3X_8,\\
(38)\;Y_3X_8=Y_4X_6-3Y_4X_4+3Y_4Y_8, & (86)\;Y_1X_4=Y_5Y_4-3Y_3Y_4+3Y_2Y_4,\\
(39)\;Y_1Y_8=Y_3Y_4-3Y_2Y_4+3Y_1X_4, & (87)\;Y_3Y_4=Y_4Y_2-3Y_3X_2+3Y_3Y_6,\\
(40)\;Y_4X_3=Y_8Y_3-3Y_7X_3+3Y_6X_3, & (88)\;Y_2X_7=Y_4X_3-3Y_4Y_7+3Y_3X_7,\\
(41)\;Y_3X_7=Y_4X_5-3Y_4X_3+3Y_4Y_7, & (89)\;Y_1X_3=Y_5Y_3-3Y_3Y_3+3Y_2Y_3,\\
(42)\;Y_1Y_7=Y_3Y_3-3Y_2Y_3+3Y_1X_3, & (90)\;Y_3Y_3=Y_4Y_1-3Y_3X_1+3Y_3Y_5,\\
(43)\;Y_4X_2=Y_8Y_2-3Y_7X_2+3Y_6X_2, & (91)\;Y_2X_6=Y_4X_2-3Y_4Y_6+3Y_3X_6,\\
(44)\;Y_1X_8=Y_2X_6-3Y_2X_4+3Y_2Y_8, & (92)\;Y_1X_2=X_5Y_2-3Y_3Y_2+3Y_2Y_2,\\
(45)\;Y_1Y_6=Y_3Y_2-3Y_2Y_2+3Y_1X_2, & (93)\;Y_1Y_4=Y_2Y_2-3Y_1X_2+3Y_1Y_6,\\
(46)\;Y_4X_1=Y_8Y_1-3Y_7X_1+3Y_6X_1, & (94)\;Y_2X_5=Y_4X_1-3Y_4Y_5+3Y_3X_5,\\
(47)\;Y_1X_7=Y_2X_5-3Y_2X_3+3Y_2Y_7, & (95)\;Y_1X_1=Y_5Y_1-3Y_3Y_1+3Y_2Y_1,\\
(48)\;Y_1Y_5=Y_3Y_1-3Y_2Y_1+3Y_1X_1, & (96)\;Y_1Y_3=Y_2Y_1-3Y_1X_1+3Y_1Y_5.
\end{array}
$$
\end{small}

\newpage
We see that $V_4(U)$ is a quadratic algebra with generators $X_1,\dots,X_8,Y_1,\dots Y_8$ and 
relations $(i),(ii)$, $(1)-(96)$. This may not be the simplest presentation of $V_4(U)$.
Observe that the generators $Y_8$ and $Y_2$ are linear combinations of other 
generators by $(i)$ and $(ii)$, so they can be removed from the generating set.

\end{section}

\bibliography{kmMay15}

\end{document}